\documentclass[11pt,a4paper]{article}
\usepackage{amsmath}
\usepackage{amssymb}
\usepackage{color}

\usepackage{graphicx}
\usepackage{amsfonts,amsmath, amssymb}

\numberwithin{equation}{section}

\usepackage{graphics}
\usepackage{epsfig}
\newtheorem{claim}{\bf \t}[part]

\newtheorem{theorem}{Theorem}[section]

\newtheorem{remark}[theorem]{Remark}

\def\v{\varepsilon}

\def\t{\theta}

\def\b{\beta}
\def\g{\gamma}
\def\d{\delta}
\def\l{\lambda}

\def\r{\rho}
\def\s{\sigma}

\def\f{\frac}
\begin{document}

	\title{ Macroscopic regularity for the Boltzmann equation}
	
	\author{ {Feimin Huang$^{a}$, 
			Yong Wang$^{a}$\footnote{Corresponding author.  \newline \indent
				Email addresses: fhuang@amt.ac.cn(Feimin Huang), yongwang@amss.ac.cn(Yong Wang)}, }
		\\
		\ \\
		{\small \it $^a$Institute of Applied Mathematics, AMSS, CAS, Beijing 100190, China}\\
	}
	
	\date{ }
	\maketitle
	
	\begin{abstract}
		
		The regularity of solutions to the Boltzmann equation is a fundamental problem in the kinetic theory. In this paper, the case with angular cut-off is investigated.  It is shown that the macroscopic parts  of  solutions to the Boltzmann equation, i.e. the density, momentum and total energy  are continuous functions  of $(x,t)$ in the region $\mathbb{R}^3\times(0,+\infty)$. More precisely,  these macroscopic quantities immediately become continuous in any positive time even though they are initially discontinuous and the discontinuities of solutions propagate only in the microscopic level. It should be noted that such kind of phenomenon  can not happen for the compressible Navier-Stokes equations in which the initial discontinuities of the density never vanish in any finite time, see \cite{Hoff}. This hints that the Boltzmann equation has better regularity effect in the macroscopic level than compressible Navier-Stokes equations.
		
		\
		
		
		\
		
		AMS: 35Q35, 35B65, 76N10
	\end{abstract}

	\section{Introduction}
	
	In this paper, we consider the Boltzmann equation
	\begin{equation}\label{1.1}
	F_t+v\cdot\nabla_x F=Q(F,F),
	\end{equation}
	where $F(t,x,v)$ is a distribution function for the gas particles at  time $t>0$,  position $x\in\mathbb{R}^3$ and particle velocity $v\in\mathbb{R}^3$. The collision operator $Q(F,F)$ takes the following bilinear form 
	\begin{align}\label{1.2}
	Q(F_1,F_2)&=\int_{\mathbb{R}^3}\int_{\mathbb{S}^2} q(\t,v-u)F_1(u')F_2(v')dud\omega-\int_{\mathbb{R}^3}\int_{\mathbb{S}^2} q(\t,v-u)F_1(u)F_2(v)dud\omega\nonumber\\
	&\triangleq Q_+(F_1,F_2)-Q_-(F_1,F_2),
	\end{align}
	where the relation between the post-collison velocity $(v',u')$ of the two particles with pre-collision velocity pair $(v,u)$ is given by
	\begin{equation*}\label{1.3}
	u'=u+[(v-u)\cdot\omega]\omega,~~v'=v-[(v-u)\cdot\omega]\omega,
	\end{equation*}
	satisfying
	$$u'+v'=u+v,~~\mbox{and}~~|u'|^2+|v'|^2=|u|^2+|v|^2.$$
	Throughout this paper,  we consider the hard and soft potentials with angular cut-off, i.e.,
	\begin{equation}\label{1.4}
	q(\t,v-u)=|v-u|^{\gamma}b(\t),~\mbox{and}~\int_0^{\f\pi2}b(\t)d\t<+\infty,~~\mbox{with}~~-3<\gamma\leq1.
	\end{equation} 
	We impose the Boltzmann equation \eqref{1.1} with initial condition
	\begin{equation}\label{1.5-1}
	F(t,x,v)|_{t=0}=F_0(x,v)\geq 0.
	\end{equation}
	
	There is an enormous literature on the  study of well-posedness of the Cauchy problem for the Boltzmann equation.  The existence of  global renormalized solutions was first proved  by R. Diperna and P.-L. Lions \cite{D-Lion} for large initial data. However, little is known on the uniqueness of renormalized solution.  For small initial data, various results on the existence and uniqueness of solutions were obtained, see  \cite{Bellomo-T,Duan2,Guo1,Guo,Ham,IIIner-S,Imai-N,Kaniel-S,Perthame,Pole,Ukai-Yang} and references therein.  We also refer  to the monographs \cite{Bellomo,Cercignani,Glassey} and  survey paper \cite{Villani}.
	
	On the other hand, the regularity of solutions is a key problem of the Boltzmann equation.  For the case with angular cut-off, it is believed that  the initial singularities propagate in time since the Boltzmann equation is hyperbolic.  This property was proved by Boudin-Desvillettes \cite{Desvillettes}  with   propagation of Sobolev $H^{\f1{25}}$ singularity  in the case near vacuum,  later by Duan-Li-Yang \cite{Duan}  in the case near global Maxwellian. We also refer to \cite{Aoki,Desvillettes1,Kim,Villani1} for other interesting works.  For the case without angular cut-off, the solution is believed  more regular  since the non-cutoff collision operator behaves like a fractional Laplacian. For interesting results on this topic, we refer to \cite{Alex-Yang,Alex-Yang1,Can-K,Carlen,Desvillettes-W,M-W-Yang} and the references therein.

	It is well known that the Boltzmann equation is closely related to the compressible Navier-Stokes equations through Chapman-Enskog expansion. It is conjectured that the Boltzmann equation has regularity in the macroscopic components. In fact, the famous averge lemma \cite{golse} shows that $\int_{R^3}f(t,x,v)\varphi(v)dv\in H^{\frac12}(t,x)$ for any smooth compact support function $\varphi(v)$. This indicates that the macroscopic components like density, momentum and total energy probably have $H^{\frac12}(t,x)$ regularity. However more regularity is not known from the average lemma. The purpose of this paper is to investigate better regularity for the macroscopic quantities of solutions to Boltzmann equation \eqref{1.1} with angular cut-off.  We find that the Boltzmann equation has better macroscopic regularity beyond the compressible Navier-Stokes equations.  Precisely speaking,  we shall prove  that the macroscopic parts  of  solutions to the Boltzmann equation, i.e. the density, momentum and total energy  are continuous functions  of $(x,t)$ in the region $\mathbb{R}^3\times(0,+\infty)$. In other words,  these macroscopic quantities immediately become continuous in any positive time even though they are initially discontinuous and the discontinuities of solutions propagate only in the microscopic level. It should be noted that such kind of phenomenon  can not happen for the compressible Navier-Stokes equations in which the initial discontinuities of the density never vanish in any finite time, see \cite{Hoff}. This means that the Boltzmann equation has better regularity effect in the macroscopic level beyond the compressible Navier-Stokes equations. 
	
	\
	
	Define a weight function 
	\begin{equation}\label{WF}
	w(v)=(1+|v|^2)^{\f\beta2},
	\end{equation}
	and let the global Maxwellian $\mu$ be 
	\begin{equation}\label{GM}
	\mu(v)=\f{1}{(2\pi)^{\f32}}\exp\left(-\f{|v|^2}{2}\right).
	\end{equation}
	Let $F_0$ be the initial data satisfying,
	{\it  for any fixed $(x,t)\in\mathbb{R}^3\times(0,+\infty)$}, 
	\begin{equation}\label{1.5}
	\|wF_0\|_{L^\infty_{x,v}}<+\infty~\mbox{and}~\int_{\mathbb{R}_{loc}^3}|F_0(x+\d-v(t+\s),v)-F_0(x-vt,v)|dv\rightarrow 0,~\mbox{as}~\s, \d\rightarrow0.
	\end{equation}

	We note that the initial condition \eqref{1.5} is general.
	Indeed, there is a large class of  initial data satisfying \eqref{1.5}. For the case near global Maxwellian, we can choose, for example, that
	\begin{align}\label{1.8}
	F_0(x,v)=\f{\r_0(x)}{(2\pi\t_0(x))^{\f32}}(1+P(v))\exp{\Big(-\f{|v-u_0(x)|^2}{2\t_0(x)}\Big)},
	\end{align}
	where $|P(v)|\leq C(1+|v|)^m$ holds for some positive constant $m>0$, which is allowed to be discontinuous in $v\in R^3$, $\r_0, \t_0, u_0$ may also be discontinuous in  $x\in R^3$ and satisfy
	\begin{align}\label{1.9}
	0\leq \r_0(x)\leq \hat{C}<+\infty,~~\hat{C}^{-1}\leq \t_0(x)\leq \hat{C}<+\infty,~~|u_0(x)|\leq \hat{C}<+\infty,
	\end{align}
	with
	\begin{align}\label{1.10}
	\|(\r_0-1,\t_0-1,u_0)\|_{L^p}<+\infty,~\mbox{for some}~p\in[1,\infty).
	\end{align}
For the case near vacuum, one can choose, for example, that
	\begin{align}\label{1.8-1}
	F_0(x,v)=\f{1}{(2\pi)^{\f32}}\r_0(x)(1+P(v))\exp{\Big(-\f{|v-u_0(x)|^2}{2}\Big)},
	\end{align}
	where  $\r_0$ and $ u_0$ satisfy
	\begin{align}\label{1.9-1}
	\begin{cases}
	0\leq \r_0(x)\leq \hat{C}<+\infty,~~|u_0(x)|\leq \hat{C}<+\infty,\\[1.5mm]
	\|(\r_0, u_0)\|_{L^p}<+\infty,~\mbox{for some}~p\in[1,\infty).
	\end{cases}
	\end{align}
One can check that  both   \eqref{1.8} and \eqref{1.8-1}  satisfy the condition \eqref{1.5}, see appendix below.
	
	\
	
The solutions considered in this paper are in the following space {\bf X}:\\[2mm]
{\it Near Maxwellian $\mu$,} $F(t,x,v)\geq0$ for a.e. $(t,x,v)\in(0,+\infty)\times\mathbb{R}^3\times\mathbb{R}^3$ and 
	\begin{equation}\label{1.6}
	w(\cdot)[F(\cdot,\cdot,\cdot)-\mu(\cdot)]\in L^{\f{b}{b-1}}(0,T;L^{2}_{x,v})~~\mbox{and}~
	w(\cdot)F(\cdot,\cdot,\cdot)\in L^b(0,T; L^\infty_{x,v})~\mbox{for}~b>2;
	\end{equation}
 {\it  Near vacuum, } $F(t,x,v)\geq0$ for a.e. $(t,x,v)\in(0,+\infty)\times\mathbb{R}^3\times\mathbb{R}^3$ and 
	\begin{equation}\label{1.6-1}
	w(\cdot)F(\cdot,\cdot,\cdot) \in  L^{\f{b}{b-1}}(0,T;L^{2}_{x,v})~~\mbox{and}~
	w(\cdot)F(\cdot,\cdot,\cdot)\in L^b(0,T; L^\infty_{x,v})~\mbox{for}~b>2.
	\end{equation}


\

Then our main results are as follows:
\begin{theorem}\label{thm1.1}
Let \eqref{1.4} hold, $\b>8+2\max\{0,\g\}$, the initial data $F_0$ satisfy \eqref{1.5},  and  the solutions $F(t,x,v)$ of  Boltzmann  equation belong to the space  {\bf X}, then the macroscopic components of  solutions to the Boltzmann equation \eqref{1.1}  are continuous functions of $(x,t)$ in the region $\mathbb{R}^3\times(0,+\infty)$.
\end{theorem}

\

Moreover, we have 
\begin{theorem}\label{thm1.2}
Under the conditions of Theorem \ref{thm1.1}, and if the initial data $F_0$ further satisfies
\begin{equation}\label{1.5-2}
\sup_{(x,t)\in\mathbb{R}^3\times[t_1,T]}\int_{\mathbb{R}^3_{loc}}|F_0(x+\d-v(t+\s),v)-F_0(x-vt,v)|dv\rightarrow 0,~\mbox{as}~\s,\d\rightarrow0.
\end{equation}
where $t_1$ and $T$ are any  given times with $0<t_1<T<+\infty$, then the macroscopic components of  solutions to the Boltzmann equation \eqref{1.1} are uniformly continuous functions of $(x,t)$ in $\mathbb{R}^3\times[t_1,T]$. 
\end{theorem}

A few remarks are in order.

\begin{remark}
It is noted that  the condition \eqref{1.5-2} is slightly stronger than \eqref{1.5}.  Indeed, both examples  in \eqref{1.8} and \eqref{1.8-1} satisfy \eqref{1.5-2} for any $0<t_1<T<+\infty$, see the appendix below.
\end{remark}

\begin{remark} 
It is noted that even though the initial macroscopic parts in Theorems \ref{thm1.1} and \ref{thm1.2} are allowed to be discontinuous in $x\in \mathbb{R}^3$ (the possible discontinuous set has zero measure),  the macroscopic components of solutions to the Boltzmann equation immediately  become  continuous in $(x,t)$ for any positive time $t>0$.  Such kind of phenomenon  can not happen for the compressible Navier-Stokes equations in which the initial discontinuities of the density never vanish in any finite time, see \cite{Hoff}. This hints that the Boltzmann equation has better regularity effect in the macroscopic level beyond the compressible Navier-Stokes equations, and  the discontinuities of solutions to the Boltzmann equation propagate only in the microscopic level.

\end{remark}


\begin{remark}
Theorems \ref{thm1.1} and \ref{thm1.2} are still valid if the solution space $L^{\f{b}{b-1}}(0,T;L^{2}_{x,v})$ is relaxed to $L^{\f{b}{b-1}}(0,T;L^{a}_{x,v})$ with $a>1$. 
\end{remark}

\begin{remark}
If $\b>8+2\max\{0,\g\}+k,~k\in\mathbb{R}_+$, let $B(v)$ be any function satisfying $|B(v)|\leq C(1+|v|)^{2+k}$,
 then, by the same arguments as in the proof of Theorems \ref{thm1.1} and \ref{thm1.2},  one can still prove that 
\begin{equation}\nonumber
\int_{\mathbb{R}^3} B(v) F(t,x,v)dv,
\end{equation} 
has the same results as in Theorems \ref{thm1.1} and \ref{thm1.2}.
\end{remark}

\begin{remark}
It is noted  that the solutions constructed in \cite{Ukai-Yang,Guo}  near Maxwellian and \cite{Bellomo-T,Ham,IIIner-S,Kaniel-S,Perthame} near vacuum belong to {\bf X}. If their initial data further satisfy \eqref{1.5}  (or \eqref{1.5-2}), then Theorem \ref{thm1.1} (or Theorem \ref{thm1.2})  holds. That is, the macroscopic quantities of their solutions are  continuous in $(x,t)$.  Thus, Theorems \ref{thm1.1} and   \ref{thm1.2}  are self-contained in this  sense.

\end{remark}

\

The proof of Theorem \ref{thm1.1} and Theorem \ref{thm1.2} is based on the following observation. Since the Boltzmann equation is essentially transport, the solution propagates along the velocity $v$ direction. Thus the integral with respect to $v$  can  be translated to the integral with respect to $x$ by the change of variables. The continuity of macroscopic components 
is then derived from the continuity of translations on $L_{x,t}^p$ space. 
The formula of collision operator plays important role in the whole proof.  

\

The rest of this paper is arranged as follows. Section 2 is devoted to the proof of Theorem \ref{thm1.1}.
The proof of Theorem \ref{thm1.2} is given in section 3. 

\

\noindent{\bf Notations.}  Throughout this paper, $C$ denotes a generic positive constant which may depend on $\b,\g$ and  vary from line to line.  $C(t)$ denotes the generic positive continuous function depending on time $t>0$ and $\b,\g$ which also may vary from line to line. $\|\cdot\|_{L^2}$ denotes the standard $L^2(\mathbb{R}^3_x\times\mathbb{R}^3_v)$-norm, and $\|\cdot\|_{L^\infty}$ denotes the $L^\infty(\mathbb{R}^3_x\times\mathbb{R}^3_v)$-norm.


\section{Proof of the Main Theorems}

Let $F(t,x,v)$ be the solution of Boltzmann equation \eqref{1.1},  we write down the Duhamel form of the solution to Boltzmann equation \eqref{1.1}.  For  $(t,x,v)\in(0,+\infty)\times\mathbb{R}^3\times\mathbb{R}^3$, one has that 
\begin{align}\label{2.1}
F(t,x,v)&=F_0(x-vt,v)e^{-\int_0^t g(\tau;t,x,v)d\tau}+\int_0^te^{-\int_s^t g(\tau;t,x,v)d\tau} Q_+(F,F)(s,x-v(t-s),v)ds\nonumber\\
&\triangleq I_1(t,x,v)+I_2(t,x,v),
\end{align}
where
\begin{equation}\label{2.2}
g(\tau;t,x,v)=\int_{\mathbb{R}^3}\int_{\mathbb{S}^2} q(\t,v-u) F(\tau,x-v(t-\tau),u)d\omega du\geq 0.
\end{equation}

\

\

\noindent{\bf Proof of Theorem \ref{thm1.1}.}  For any fixed $(t,x)\in(0,+\infty)\times\mathbb{R}^3$, let $\s\in\mathbb{R}$, $\d\in\mathbb{R}^3$ with $|\s|\leq \min\{1,\f{t}{2}\}$ and $|\d|\leq 1$. The continuity of density $\rho(x,t)$ in $(x,t)$ is equivalent to show that 
\begin{align}\label{2.3}
\Big|\int_{\mathbb{R}^3}F(t+\s,x+\d,v)-F(t,x,v)dv\Big|\rightarrow0~\mbox{as}~\s, \d\rightarrow0.
\end{align}
That is, for any $\v>0$, we need to prove that there exists $\chi>0$, which may depend on $(x,t)$ and $\v$, such that if $|\s|\leq \chi$ and $|\d|\leq \chi$, 
\begin{align}\label{2.3-1}
\Big|\int_{\mathbb{R}^3}F(t+\s,x+\d,v)-F(t,x,v)dv\Big|\leq \v.
\end{align}
We assume $\s\geq0$ without loss of generality and denote $t_\s\doteq t+\s$ for notation simplicity.  We divide the proof into two parts. The proof is based on 
the initial condition \eqref{1.5} and  the continuity of translations on $L^p$.\\

\noindent{\bf Part 1: Estimation on $\int_{\mathbb{R}^3}I_1(t,x,v)dv$.}\\
Firstly, we note that 
\begin{align}\label{2.4}
&J_1\triangleq\Big|\int_{\mathbb{R}^3}(I_1(t_\s,x+\d,v)-I_1(t,x,v))dv\Big|\nonumber\\
&\leq \int_{\mathbb{R}^3}\Big|F_0(x+\d-vt_\s,v)-F_0(x-vt,v)\Big|dv+\int_{\mathbb{R}^3}F_0(x-vt,v)\Big|\int_t^{t_\s}g(\tau;t_\s,x+\d,v)d\tau\Big|dv\nonumber\\
&~~~~~~~+\int_0^t \int_{\mathbb{R}^3}F_0(x-vt,v) \Big|g(\tau;t_\s,x+\d,v)-g(\tau;t,x,v)\Big| dvd\tau\nonumber\\
&\triangleq J_{11}+J_{12}+J_{13}.
\end{align}
Let $\tilde{F}(\s,\d)=F_0(x+\d-vt_\s,v)-F_0(x-vt,v)$, it follows from  \eqref{1.5} that 
\begin{align}\label{2.5}
J_{11}
&\leq \int_{|v|\geq N}\Big|\tilde{F}(\s,\d)\Big|dv+\int_{|v|\leq N}\Big|\tilde{F}(\s,\d))\Big|dv\nonumber\\
&\leq C\|w F_0\|_{L^\infty}\int_{|v|\geq N}(1+|v|)^{-\b}dv+\int_{|v|\leq N}\Big|\tilde{F}(\s,\d)\Big|dv\nonumber\\
&\leq C(1+|N|)^{-\b+3}+\int_{|v|\leq N}\Big|\tilde{F}(\s,\d)\Big|dv.
\end{align}
Fix  $N$ large such that $
C(1+|N|)^{-\b+3}\leq \f\v{32},~\mbox{for}~\b>3.
$
Then, from  \eqref{1.5},  there exists $\chi_1>0$ depending on $(x,t)$ and $\v$ such that if $\s+|\d|\leq \chi_1$, it holds that 
\begin{align}\label{2.5-3}
J_{11}\leq \f{\v}{16}.
\end{align}
For $J_{12}$, it follows from \eqref{2.2} that  
\begin{align}\label{2.5-4}
J_{12}&\leq C\|w F_0\|_{L^\infty}\int_{\mathbb{R}^3}w(v)^{-1}\int_t^{t+\s}\int_{\mathbb{R}^3}|v-u|^{\g}w(u)^{-1}\|w F(\tau)\|_{L^\infty}dud\tau dv\nonumber\\
&\leq C\int_t^{t+\s}\|w F(\tau)\|_{L^\infty}d\tau\int_{\mathbb{R}^3}\int_{\mathbb{R}^3}|v-u|^{\g}w(v)^{-1}w(u)^{-1}dudv\nonumber\\
&\leq C\s^{1-\f1b}\Big(\int_t^{t+\s}\|w F(\tau)\|^b_{L^\infty}d\tau\Big)^{\f1b}\leq C(t)\s^{1-\f1b}\rightarrow0~\mbox{as}~\s\rightarrow0,
\end{align}
where we have used the fact that $F\in {\bf X}$ and 
\begin{align}\label{2.8}
\int_{\mathbb{R}^3}\int_{\mathbb{R}^3} w(v)^{-1}w(u)^{-1} |v-u|^{\gamma}dudv\leq C<+\infty,~~\mbox{for}~\beta>3+\max\{0,\g\}.
\end{align}
Hence there exists $\chi_2>0$ depending  on $t, \v$ such that if $\s\leq \chi_2$, it holds that 
\begin{equation}\label{2.5-5}
J_{12}\leq \f{\v}{16}.
\end{equation}
For $J_{13}$,  we note that 
\begin{align}\label{2.6-1}
 J_{13}&= \int_0^t \int_{\mathbb{R}^3}F_0(x-vt,v) \Big|g(\tau;t_\s,x+\d,v)-g(\tau;t_\s,x,v)\Big| dvd\tau\nonumber\\
 &~~~+\int_0^t \int_{\mathbb{R}^3}F_0(x-vt,v) \Big|g(\tau;t_\s,x,v)-g(\tau;t,x,v)\Big| dvd\tau\nonumber\\
 &\triangleq J_{131}+J_{132}.
\end{align}
It holds that 
\begin{align}\label{2.6}
 J_{131}
& \leq  C\int_0^t \int_{\mathbb{R}^3}\int_{\mathbb{R}^3} w(v)^{-1}w(u)^{-1} |v-u|^{\gamma}\nonumber\\ 
&~~~~~~~~~\times \Big|w(u)[F(\tau,x+\d-v(t_\s-\tau),u)-F(\tau,x-v(t_\s-\tau),u)]\Big|du dvd\tau\nonumber\\
&=\Big\{\int_{t-\l}^{t}+\int_0^{t-\l}\Big\}(\cdots)d\tau\triangleq J_{1311}+J_{1312}.
\end{align}
From \eqref{2.8}, one has that 
\begin{align}\label{2.7}
J_{1311} \leq C \int_{t-\l}^{t}\|wF(\tau)\|_{L^\infty}d\tau\leq C(t)\l^{1-\f1b}.
\end{align}
On the other hand, it holds for some $p>1$ that 
\begin{align}\label{2.9}
&J_{1312}
\leq C\int_0^{t-\l} \Big(\int\int \Big|w(u)[F(\tau,x+\d-v(t_\s-\tau),u)-F(\tau,x-v(t_\s-\tau),u)]\Big|^{p_\ast}du dv\Big)^{\f1{p_\ast}}\nonumber\\ 
&~~~~~~~~~~~~~~~~~~~~~~~\times \Big(\int\int w(v)^{-p}w(u)^{-p} |v-u|^{p\gamma}dudv\Big)^{\f1p} d\tau,
\end{align}
where $p_\ast=\f{p}{p-1}$. Since $|v-u|^{p\g}$ is singular when $-3<\g<0$, we choose 
\begin{eqnarray}\label{2.9-2}
p=1+\f{3+\g}{4(9-\g)},~~\mbox{for}~~-3<\g\leq 1,
\end{eqnarray}  
so that $
1<p\leq \f87,~~p_\ast=\f{p}{p-1}\geq7,~~p\g>-3,~~p\f{3-\g}{2}>3$
and 
\begin{align}\label{2.9-3}
\int_{\mathbb{R}^3}\int_{\mathbb{R}^3}|v-u|^{p\g}w(u)^{-p}w(v)^{-p}dudv\leq C<\infty~~\mbox{for}~~\b>3+\max\{0,\g\}.
\end{align}
Changing variable by $y=x-v(t_\s-\tau)$, we obtain, for any fixed $\l>0$,  that
\begin{align}\label{2.9-4}
&J_{1312}
\leq C \l^{-\f3{p_\ast}}\int_0^{t-\l}\|w F(\tau)\|^{1-\f2{p_\ast}}_{L^\infty}\|w(\cdot)[F(\tau,\cdot+\d,\cdot)-F(\tau,\cdot,\cdot)]\|_{L^2}^{\f2{p_\ast}} d\tau\nonumber\\
&\leq C \l^{-\f3p_\ast}\Big(\int_0^{t-\l}\|\cdot\|^b_{L^\infty}d\tau\Big)^{\f1b-\f{2}{bp_{\ast}}}
\Big(\int_0^{t-\l}\|\cdot\|_{L^2}^{\f{b}{b-1}}d\tau\Big)^{\f{2(b-1)}{bp_{\ast}}}\Big(\int_0^{t-\l}d\tau\Big)^{\tilde{d}}
\nonumber\\
&\leq C(t)\l^{-\f3p_\ast}\Big(\int_0^{t-\l}\|w(\cdot)[F(\tau,\cdot+\d,\cdot)-F(\tau,\cdot,\cdot)]\|_{L^2}^{\f{b}{b-1}}d\tau\Big)^{\f{2(b-1)}{bp_{\ast}}}\rightarrow0,
~\mbox{as}~\d\rightarrow0,
\end{align}
where  $\tilde{d}={\small {\tiny }}1-\f1b+\f{4}{bp_{\ast}}-\f2{p_\ast}>0$ due to $b>2,~p_\ast\geq7$.
From \eqref{2.6}, \eqref{2.7} and \eqref{2.9-4},  choosing $\l^{1-\f1b}= \f{\v}{64C(t)}$, then there exists $\chi_3>0$ depending only on $t, \v$ such that  if $\s+|\d|\leq \chi_3$,  one obtains, for $\b>3+\max\{0,\g\}$, that 
\begin{equation}\label{2.10}
J_{131}\leq \f\v{16}.
\end{equation}

For $J_{132}$, one notes that 
\begin{align}\label{3.16}
 J_{132}
&\leq C\|wF_0\|_{L^\infty}\int_0^t \int_{\mathbb{R}^3} \int_{\mathbb{R}^3}w(v)^{-1}w(u)^{-1}|v-u|^{\gamma}\nonumber\\
&~~~~~~~~~~~~~~~~\times \Big| w(u)[F(\tau,x-v(t_\s-\tau),u)- F(\tau,x-v(t-\tau),u)]\Big|dudvd\tau\nonumber\\
&\leq C\Big\{\int_{0}^t\int_{|v|\geq N} \int_{\mathbb{R}^3}+\int_{t-\l}^t\int_{|v|\leq N} \int_{\mathbb{R}^3}+\int_0^{t-\l}\int_{|v|\leq N} \int_{\mathbb{R}^3}\Big\}(\cdots)dvdud\tau\nonumber\\
&\triangleq J_{1321}+J_{1322}+J_{1323}.
\end{align}
It is straightforward to obtain, for $\b>3+\max\{0,\g\}$,  that 
\begin{align}\label{3.23}
J_{1321}\leq C\int_0^t\int_{|v|\geq N}w(v)^{-1}(1+|v|)^{\max\{0,\g\}}dvd\tau
\leq C(t)(1+N)^{-\b+3+\max\{0,\g\}}\leq \f{\v}{64},
\end{align}
and
\begin{align}\label{3.17}
J_{1322}\leq C\int_{t-\l}^t\|wF(\tau)\|_{L^\infty}d\tau\leq C(t)\l^{1-\f1b}\leq \f{\v}{64},
\end{align}
where $N$ and $\l$ are suitably chosen so that $C(t)(1+N)^{-\b+3+\max\{0,\g\}}\leq \f{\v}{64}$ and $\l^{1-\f1b}=\f{\v}{64 C(t)}$.
Now we estimate $J_{1323}$. By the definition of {\bf X}, there exists a smooth compact support function $F^\v(\cdot,\cdot,\cdot)$ such that 
\begin{align}\label{3.18}
&\left(\int_0^t\|w(F(\tau,\cdot,\cdot)-F^\v(\tau,\cdot,\cdot))\|_{L^2}^{\f{b}{b-1}}d\tau\right)^{\f{2(b-1)}{bp_\ast}}
\leq \f{\l^{\frac{3}{p_*}}\v}{C_1},
\end{align}
where $C_1$ will be chosen later.
Using \eqref{3.18} and the same arguments as in \eqref{2.9-4}, for $\b>3+\max\{0,\g\}$, we have 
\begin{align}\label{3.19}
	&J_{1323}\leq C\int_0^{t-\l} \int_{|v|\leq N} \int_{\mathbb{R}^3}w(v)^{-1}w(u)^{-1}|v-u|^{\g}\nonumber\\
	&~~~~~~~~~~~~\times  \Big\{\Big|w(u)[F(\tau,x-v(t_\s-\tau),u)- F^\v(\tau,x-v(t_\s-\tau),u)]\Big|\nonumber\\
	&~~~~~~~~~~~~+ \Big|w(u)[F(\tau,x-v(t-\tau),u)- F^\v(\tau,x-v(t-\tau),u)]\Big|\nonumber\\
	&~~~~~~~~~~~~+ \Big|w(u)[F^\v(\tau,x-v(t_\s-\tau),u)- F^\v(\tau,x-v(t-\tau),u)]\Big|\Big\}dudvd\tau\nonumber\\
	&\leq C(t)\l^{-\f3p_\ast}\left(\int_0^t\|w(F(\tau,\cdot,\cdot)-F^\v(\tau,\cdot,\cdot))\|^{\f{b}{b-1}}_{L^2}d\tau\right)^{\f{2(b-1)}{bp_\ast}}\nonumber\\
	&~~~~~~~+C_{N,\v} t\sup_{|v|\leq N,(\tau,y,u)\in((0,T]\times\mathbb{R}^3\times\mathbb{R}^3)}\Big|w(u)[F^\v(\tau,y-v\s,u)-F^\v(\tau,y,u)]\Big|\nonumber\\
	&\leq \f{C(t)}{C_1}\v+\f{\v}{128}\leq \f{\v}{64},
	\end{align}
where we have  used the fact that
\begin{equation}\label{3.21}
\lim_{\s\rightarrow0}C_{N,\v} t\sup_{|v|\leq N,(\tau,y,u)\in((0,t]\times\mathbb{R}^3\times\mathbb{R}^3)}\Big|w(u)[F^\v(\tau,y-v\s,u)-F^\v(\tau,y,u)]\Big|=0,
\end{equation}
and $C_1$ is chosen large so that \eqref{3.19} holds. Thus, from \eqref{3.16}-\eqref{3.19},  there exists  $\chi_4>0$ depending on $t$ and $\v$ such that if $\s\leq \chi_4$, it holds that 
\begin{equation}\label{3.22}
J_{132}\leq \f{\v}{16}.
\end{equation}
Therefore, for $\b>3+\max\{0,\g\}$, it follows from \eqref{2.4}, \eqref{2.5-3}, \eqref{2.5-5}, \eqref{2.6-1}, \eqref{2.10} and \eqref{3.22}  that 
\begin{align}\label{3.22-1}
\Big|\int_{\mathbb{R}^3}\Big\{I_1(t+\s,x+\d,v)-I_1(t,x,v)\Big\}dv\Big|\leq \f{\v}{2},~~\mbox{for}~\s+|\d|\leq \min\{\chi_1,\chi_2,\chi_3,\chi_4\}.
\end{align}

\

\noindent{\bf Part 2: Estimation on $\int_{\mathbb{R}^3}I_2(t,x,v)dv$.}\\
Note  that 
\begin{align}\label{2.11-1}
 J_2&\triangleq \Big|\int_{\mathbb{R}^3}I_2(t_\s,x+\d,v)-I_2(t,x,v)dv\Big|\nonumber\\
&\leq \Big|\int_{\mathbb{R}^3}I_2(t_\s,x+\d,v)-I_2(t_\s,x,v)dv\Big|+\Big|\int_{\mathbb{R}^3}I_2(t_\s,x,v)-I_2(t,x,v)dv\Big|\nonumber\\
&\triangleq J_{21}+J_{22}.
\end{align}
For $J_{21}$, it follows from \eqref{2.1}  that
\begin{align}\label{2.11}
 J_{21}&\leq \int_0^{t_\s}\int_{\mathbb{R}^3}Q_+(F,F)(s,x+\d-v(t_\s-s),v)\int_s^{t_\s}|g(\tau;t_\s,x+\d,v)-g(\tau;t_\s,x,v))|d\tau dvds\nonumber\\
&~~~+\int_0^{t_\s}\int_{\mathbb{R}^3}\Big|Q_+(F,F)(s,x+\d-v(t_\s-s),v)-Q_+(F,F)(s,x-v(t_\s-s),v)|dvds\nonumber\\
&\triangleq J_{211}+J_{212}.
\end{align}
A direct calculation shows  that 
\begin{align}\label{2.12}
Q_+(F,F)(s,y,v)
&\leq \int_{\mathbb{R}^3}\int_{\mathbb{S}^2}q(\t,v-u)w(u')^{-1}w(v')^{-1}|w(u')F(s,y,u')|\cdot|w(v')F(s,y,v')|dud\omega\nonumber\\
&\leq  \|wF(s)\|^2_{L^\infty}\int_{\mathbb{R}^3}(1+|v|+|u|)^{-\beta}|v-u|^{\gamma}du\nonumber\\
&\leq C(1+|v|)^{-\b+3+\max\{0,\g\}}\|wF(s)\|^2_{L^\infty},~~\mbox{for}~~\b>3+\max\{0,\g\},
\end{align}
where we have used the fact that 
\begin{equation}\label{2.12-1}
w(u')w(v')\geq (1+|v|^2+|u|^2)^{\f\beta2},~\b\geq0,
\end{equation}
and
\begin{align}\label{2.12-3}
\int_{\mathbb{R}^3}(1+|v|+|u|)^{-\beta}|v-u|^{\gamma}du 
\leq C(1+|v|)^{-\b+3+\max\{0,\g\}}~~\mbox{for}~\b>3+\max\{0,\g\}.
\end{align}
Thus, it follows from \eqref{2.2} and \eqref{2.12} that 
\begin{align}\label{2.13}
 J_{211}&\leq C\int_0^{t_\s} \|wF(s)\|^2_{L^\infty}ds\cdot\int_0^{t_\s}\int_{\mathbb{R}^3}\int_{\mathbb{R}^3}|v-u|^\g w(u)^{-1}\cdot(1+|v|)^{-\b+3+\max\{0,\g\}} \nonumber\\
&~~~~~~~~~~~~~~~~~\times \Big|w(u)[F(\tau,x+\d-v(t_\s-\tau),u)-F(\tau,x-v(t_\s-\tau),u)]\Big|du dv d\tau \nonumber\\
&\leq C(t)\int_0^{t_\s}\int_{\mathbb{R}^3}\int_{\mathbb{R}^3}|v-u|^\g w(u)^{-1}\cdot(1+|v|)^{-\b+3+\max\{0,\g\}} \nonumber\\
&~~~~~~~~~~~~~~~~~\times \Big|w(u)[F(\tau,x+\d-v(t_\s-\tau),u)-F(\tau,x-v(t_\s-\tau),u)]\Big|du dv d\tau \nonumber\\
&\leq C(t)\cdot\Big\{\int_{t-\l}^{t_\s}+\int_0^{t-\l}\Big\}(\cdots)d\tau\triangleq J_{2111}+J_{2112}.
\end{align}
It is straightforward to obtain that 
\begin{align}\label{2.13-1}
J_{2111}&\leq C(t) \int_{t-\l}^{t_\s}\|w F(\tau)\|_{L^\infty}d\tau\int_{\mathbb{R}^3}(1+|v|)^{-\b+3+2\max\{0,\g\}} dv\nonumber\\
&\leq C(t)(\s+\l)^{1-\f1b},~~\mbox{for}~~\b>6+2\max\{0,\g\}.
\end{align}
On the other hand, for any fixed $\l>0$, it holds that 
\begin{align}\label{2.13-2}
J_{2112}
&\leq C(t)\int_0^{t-\l} \Big(\int_{\mathbb{R}^3}\int_{\mathbb{R}^3}  |v-u|^{p\gamma}w(u)^{-p}(1+|v|)^{-p\b+3p+p\max\{0,\g\}}dudv\Big)^{\f1p}\nonumber\\ 
&~~~~~~~\times\Big(\int_{\mathbb{R}^3}\int_{\mathbb{R}^3} \Big|w(u)[F(\tau,x+\d-v(t_\s-\tau),u)-F(\tau,x-v(t_\s-\tau),u)]\Big|^{p_\ast}du dv\Big)^{\f1{p_\ast}}d\tau\nonumber\\
&\leq C(t) \l^{-\f3p_\ast}\int_0^{t-\l}\Big(\int_{\mathbb{R}^3}\int_{\mathbb{R}^3} w(u)|F(\tau,y+\d,u)-F(\tau,y,u)|^{p_\ast}du dy\Big)^{\f1p_\ast}d\tau\nonumber\\
&\leq C(t)\l^{-\f3p_\ast}\Big(\int_0^{t-\l}\|w(u)[F(\tau,\cdot+\d,u)-F(\tau,\cdot,u)]\|_{L^2}^{\f{b}{b-1}}d\tau\Big)^{\f{2(b-1)}{bp_{\ast}}}\rightarrow0
~\mbox{as}~\d\rightarrow0,
\end{align}
where $p$ is defind in \eqref{2.9-2}. Therefore, for $\b>6+2\max\{0,\g\}$, combining  \eqref{2.13}-\eqref{2.13-2}, there exists $\chi_5>0$ depending on $t$ and $\v$ such that if $\s+|\d|\leq \chi_5$, it holds that  
\begin{equation}\label{2.14}
 J_{211}\leq \f{\v}{32}.
\end{equation}

For $J_{212}$,  we have  
\begin{align}\label{2.15}
J_{212}&\leq \int_0^{t_\s}\int_{\mathbb{R}^3}\int_{\mathbb{R}^3}\int_{\mathbb{S}^2}q(\t,v-u)F(s,x+\d-v(t_\s-s),u')\nonumber\\
&~~~~~~~~~~~\times \Big|F(s,x+\d-v(t_\s-s),v')-F(s,x-v(t_\s-s),v')\Big|dudvd\omega ds\nonumber\\
&~~+\int_0^{t_\s}\int_{\mathbb{R}^3}\int_{\mathbb{R}^3}\int_{\mathbb{S}^2}q(\t,v-u)F(s,x-v(t_\s-s),v')\nonumber\\
&~~~~~~~~~~~\times \Big|F(s,x+\d-v(t_\s-s),u')-F(s,x-v(t_\s-s),u')\Big|dudvd\omega ds\nonumber\\
&\triangleq J_{2121}+J_{2122}.
\end{align}
We only consider $J_{2121} $ since $J_{2122}$ can be treated similaryly. It follows from \eqref{2.12-1} that 
\begin{align}\label{2.17}
J_{2121}&\le\int_0^{t_\s}\int_{\mathbb{R}^3}\int_{\mathbb{R}^3}\int_{\mathbb{S}^2}q(\t,v-u)w(u')^{-1}w(v')^{-1}\|wF(s)\|_{L^\infty}\nonumber\\
&~~~~~~~~\times \Big|w(v')[F(s,x+\d-v(t_\s-s),v')-F(s,x-v(t_\s-s),v')]\Big|dudvd\omega ds\nonumber\\
&\leq C\int_{t-\l}^{t_\s}\|wF(s)\|^2_{L^\infty}\int_{\mathbb{R}^3}\int_{\mathbb{R}^3}|v-u|^{\g}
(1+|v|+|u|)^{-\b}dudv ds\nonumber\\
&~~~~+C\int_0^{t-\l}\|wF(s)\|_{L^\infty}\int_{\mathbb{R}^3}\int_{\mathbb{R}^3}\int_{\mathbb{S}^2}q(\t,v-u)w(v)^{-\f12}w(u')^{-\f12}w(v')^{-\f12}\nonumber\\
&~~~~~~~\times w(v')\Big|F(s,x+\d-v(t_\s-s),v')-F(s,x-v(t_\s-s),v')\Big|dudvd\omega ds\nonumber\\
&\triangleq J_{21211}+J_{21212}.
\end{align}
It follows from \eqref{2.12-3} that
\begin{align}\label{2.18}
J_{21211}\leq C(t)(\s+\l)^{1-\f2b}.
\end{align}
For $J_{21212}$, we follow the idea of \cite{Bellomo}. Denote
\begin{align}\label{2.16}
z=u-v,~~z_{\shortparallel}=[z\cdot\omega]\omega,~z_{\perp}=z-z_{\shortparallel},~\eta=v+z_{\shortparallel},
\end{align}
one has $
u'=v+z_{\perp},~v'=v+z_{\shortparallel}.
$
A simple calculation yields that 
\begin{equation}\label{2.16-2}
q(\t,z)\leq C|z_{\shortparallel}|\cdot(|z_{\shortparallel}|+ |z_{\perp}|)^{\gamma-1}.
\end{equation}
Then it follows from \eqref{2.12-1} and  \eqref{2.16-2}  that 
\begin{align}\label{2.22}
&J_{21212}\leq C \int_0^{t-\l}\|wF(s)\|_{L^\infty}\int_{\mathbb{R}^3}\int_{\mathbb{R}^3}\int_{\mathbb{S}^2}w(v)^{-\f12}\cdot|z_{\shortparallel}|\cdot(|z_{\shortparallel}|+ |z_{\perp}|)^{\gamma-1}\cdot w(v+z_{\perp})^{-\f12}w(v+z_{\shortparallel})^{-\f12}\nonumber\\
&~~~~~~~\times \Big|w(v+z_{\shortparallel})[F(s,x+\d-v(t_\s-s),v+z_{\shortparallel})-F(s,x-v(t_\s-s),v+z_{\shortparallel})]\Big|dz_{\perp}d|z_{\shortparallel}| d\omega dvds\nonumber\\
&~~~~\mbox{Noting}~|z_{\shortparallel}|^2 dz_{\perp}d|z_{\shortparallel}| d\omega=dz_{\perp}d\eta,\nonumber\\
&\leq C \int_0^{t-\l}\|wF(s)\|_{L^\infty}\int_{\mathbb{R}^3}\int_{\mathbb{R}^3}\int_{z_{\perp}}w(v+z_{\perp})^{-\f12}(|\eta-v|+|z_{\perp}|)^{\gamma-1}w(v)^{-\f12}w(\eta)^{-\f12}|\eta-v|^{-1}\nonumber\\
&~~~~~~~~~~~~\times \Big|w(\eta)[F(s,x+\d-v(t_\s-s),\eta)-F(s,x-v(t_\s-s),\eta)]\Big| dz_{\perp}d\eta dvds.
\end{align}
We note, for $\b>4$, that
\begin{align}\label{2.22-1}
&\int_{z_{\perp}}w(v+z_{\perp})^{-\f12}(|\eta-v|+|z_{\perp}|)^{\gamma-1}dz_{\perp}\leq |\eta-v|^{\f{\g-1}{2}}\int_{z_{\perp}}w(v+z_{\perp})^{-\f12}|z_{\perp}|^{\f{\gamma-1}{2}}dz_{\perp}\nonumber\\
&\leq |\eta-v|^{\f{\g-1}{2}}\Big\{\int_{|z_{\perp}|\leq1}|z_{\perp}|^{\f{\gamma-1}{2}}dz_{\perp}+
\int_{|z_{\perp}|\geq1}(1+|v+z_{\perp}|)^{-\f{\b}{2}}dz_{\perp}\Big\}\leq C |\eta-v|^{\f{\g-1}{2}},
\end{align}
which, together with  \eqref{2.9-2},  yields, for $\b>6$ and any fixed $\l>0$,  that 
\begin{align}\label{2.34}
J_{21212}&\leq C \int_0^{t-\l}\|wF(s)\|_{L^\infty}\int_{\mathbb{R}^3}\int_{\mathbb{R}^3}w(v)^{-\f12}w(\eta)^{-\f12}|\eta-v|^{\f{\g-3}{2}}\nonumber\\
&~~~~~~~~~~~~~~~~~~~~\times \Big|w(\eta)[F(s,x+\d-v(t_\s-s),\eta)-F(s,x-v(t_\s-s),\eta)]\Big| d\eta dvds\nonumber\\
&\leq C\l^{-\f3p_\ast}\int_0^{t}\|wF(s)\|_{L^\infty}^{2-\f2p_{\ast}}
\|w(F(s,\cdot+\d,\cdot)-F(s,\cdot,\cdot))\|_{L^2}^{\f2p_{\ast}}\nonumber\\
&~~~~~~~~~~~~~~~~~~~~~~~~~~~~~~~~~~\times
\Big(\int_{\mathbb{R}^3}\int_{\mathbb{R}^3}w(v)^{-\f{p}2}w(\eta)^{-\f{p}2}|\eta-v|^{p\f{\g-3}{2}}d\eta dv\Big)^{\f1q}ds
\nonumber\\
&\leq C(t)\l^{-\f3p_\ast}\Big(\int_0^t \|w(F(s,\cdot+\d,\cdot)-F(s,\cdot,\cdot))\|_{L^2}^{\f{b}{b-1}}ds\Big)^{\f{2(b-1)}{bp_\ast}} \rightarrow0,~\mbox{as}~\d\rightarrow0.
\end{align}
Combining \eqref{2.15}-\eqref{2.18} and \eqref{2.34}, for $\b>6+\max\{0,\g\}$, there exists  $\chi_6>0$ depending on $t$ and $\v$ such that if $\s+|\d|\leq \chi_6$, one has that 
\begin{equation}\label{2.23}
J_{212}\leq \f{\v}{32}.
\end{equation}

It remains to estimate $J_{22}$.  A direct calculation shows  that
\begin{align}\label{3.24}
J_{22}&\leq \int_t^{t+\s}\int_{\mathbb{R}^3}Q_+(F,F)(s,x-v(t_\s-s),v)dvds\nonumber\\
&~~~+\int_0^t\int_{\mathbb{R}^3}Q_+(F,F)(s,x-v(t-s),v)\int_t^{t+\s}g(\tau;t_\s,x,v)d\tau dvds\nonumber\\
&~~~+\int_0^t\int_{\mathbb{R}^3}Q_+(F,F)(s,x-v(t-s),v)\int_s^{t}\Big|g(\tau;t_\s,x,v)-g(\tau;t,x,v)\Big|d\tau dvds\nonumber\\
&~~~+\int_0^t\int_{\mathbb{R}^3}\Big|Q_+(F,F)(s,x-v(t_\s-s),v)-Q_+(F,F)(s,x-v(t-s),v)\Big|dvds\nonumber\\
&\triangleq J_{221}+J_{222}+J_{232}+J_{224}.
\end{align}
It follows from \eqref{2.12}, for $\b>6+\max\{0,\g\}$,  that 
\begin{align}\label{3.25}
J_{221} 
\leq 
C \int_t^{t+\s}\|wF(s)\|_{L^\infty}^2ds\leq C(t)\s^{1-\f2b}\rightarrow0,~~\mbox{as}~\s\rightarrow0,
\end{align}
and
\begin{align}\label{3.26}
J_{222}&\leq \int_0^{t}\int(1+|v|)^{-\b+3+\max\{0,\g\}}\|wF(s)\|_{L^\infty}^2\int_t^{t+\s}\int
|v-u|^{\g}w(u)^{-1}\|wF(\tau)\|_{L^\infty}dud\tau dvds\nonumber\\
&\leq 
C(t) \int_t^{t+\s}\|wF(s)\|_{L^\infty}ds\leq C(t)\s^{1-\f1b}\rightarrow0,~~\mbox{as}~\s\rightarrow0.
\end{align}
Using the same arguments as in \eqref{3.23}-\eqref{3.19}, we have that if $\s$ is small, 
\begin{align}\label{3.27}
J_{223}&\leq C\int_0^{t}\|wF(s)\|_{L^\infty}^2ds \int_0^{t}\int_{\mathbb{R}^3}\int_{\mathbb{R}^3}
|v-u|^{\g}w(u)^{-1}(1+|v|)^{-\b+3+\max\{0,\g\}}\nonumber\\
&~~~~~~~~~~~~~~\times\Big|w(u)[F(\tau,x-v(t+\s-\tau),u)-F(\tau,x-v(t-\tau),u)]\Big|du dvd\tau\nonumber\\
&\leq 	C(t)\int_{t-\l}^t\|wF(\tau)\|_{L^\infty}d\tau+C(t) \int_{|v|\geq N}\int_{\mathbb{R}^3}	|v-u|^{\g}w(u)^{-1}(1+|v|)^{-\b+3+\max\{0,\g\}}dudv\nonumber\\
&~~~+C(t)\int_{0}^{t-\l}\int_{|v|\leq N}\int_{\mathbb{R}^3_u} |v-u|^{\g}w(u)^{-1}(1+|v|)^{-\b+3+\max\{0,\g\}}\nonumber\\
&~~~~\times w(u)\Big|F(\tau,x-v(t+\s-\tau),u)-F(\tau,x-v(t-\tau),u)\Big|du dvd\tau\leq \f{\v}{16}.
\end{align}
For $J_{224}$, one notes that
\begin{align}\label{3.28}
J_{224}&\leq \int_0^t\int_{\mathbb{R}^3}\int_{\mathbb{R}^3}\int_{\mathbb{S}^2}q(\t,v-u)F(s,x-v(t+\s-s),u')\nonumber\\
&~~~~~~~~~~~~~~~~\times |F(s,x-v(t+\s-s),v')-F(s,x-v(t-s),v')|dudvd\omega ds\nonumber\\
&~~+\int_0^t\int_{\mathbb{R}^3}\int_{\mathbb{R}^3}\int_{\mathbb{S}^2}q(\t,v-u)F(s,x-v(t-s),v')\nonumber\\
&~~~~~~~~~~~~~~~~\times |F(s,x-v(t+\s-s),u')-F(s,x-v(t-s),u')|dudvd\omega ds\nonumber\\
&\triangleq J_{2241}+J_{2242}.
\end{align}
We only consider $J_{2241}$ since the argument for $J_{2242}$ is similar.  It holds that
\begin{align}\label{3.29}
J_{2241}
&\leq C\int_{t-\l}^t\|wF(s)\|^2_{L^\infty}\int_{\mathbb{R}^3}\int_{\mathbb{R}^3}|v-u|^{\g}(1+|v|+|u|)^{-\b}dudvds\nonumber\\
&~~~~+\int_0^{t-\l}\|wF(s)\|^2_{L^\infty}ds\int_{|v|\geq N}\int_{\mathbb{R}^3}|v-u|^{\g}(1+|v|+|u|)^{-\b}dudv\nonumber\\
&~~~~+\int_0^{t-\l}\|wF(s)\|_{L^\infty}\int_{|v|\leq N}\int_{\mathbb{R}^3}\int_{\mathbb{S}^2}q(\t,v-u)w(u')^{-1}w(v')^{-1}\nonumber\\
&~~~~~~~~~~~\times \Big|w(v')[F(s,x-v(t+\s-s),v')-F(s,x-v(t-s),v')]\Big|dudvd\omega ds\nonumber\\
&\triangleq J_{22411}+J_{22412}+J_{22413}.
\end{align}
By the same arguments as in \eqref{3.23} and \eqref{3.17}, one can get that 
\begin{align}\label{3.30}
J_{22411}+J_{22412}\leq C(t)\l^{1-\f2b}+C(t)(1+N)^{-\b+6+\max\{0,\g\}}\leq \f{\v}{16}.
\end{align}
For $J_{22413}$, we notice that 
\begin{align}\label{3.31}
J_{22413}&\leq \int_0^{t-\l}\|wF(s)\|_{L^\infty}\int_{|v|\leq N}\int_{\mathbb{R}^3}\int_{\mathbb{S}^2}q(\t,v-u)w(v)^{-\f12}w(u')^{-\f12}w(v')^{-\f12}\nonumber\\
&~~~~~~~~\times \Big\{\Big|w(v')[F(s,x-v(t+\s-s),v')-F^\v(s,x-v(t+\s-s),v')]\Big|\nonumber\\
&~~~~~~~~~~~~~+\Big|w(v')[F(s,x-v(t-s),v')-F^\v(s,x-v(t-s),v')]\Big|\Big\}dudvd\omega ds\nonumber\\
&~~~+\int_0^{t-\l}\|wF(s)\|_{L^\infty}\int_{|v|\leq N}\int_{\mathbb{R}^3}\int_{\mathbb{S}^2}q(\t,v-u)w(u')^{-1}w(v')^{-1}\nonumber\\
&~~~~~~~~\times \Big|w(v')[F^\v(s,x-v(t+\s-s),v')-F^\v(s,x-v(t-s),v')]\Big|dudvd\omega ds\nonumber\\
&\triangleq J_{224131}+J_{224132}.
\end{align}
It follows from \eqref{3.18} and the same arguments as in \eqref{2.22}-\eqref{2.34}, for $\b>6$, that 
\begin{align}\label{3.32}
J_{224131}&\leq \l^{-\f3p_\ast}\int_0^t\|wF(s)\|_{L^\infty}^{2-\f2{p_\ast}}\|w[F(\tau,\cdot,\cdot)-F^\v(\tau,\cdot,\cdot)]\|^{\f2{p_\ast}}_{L^2}d\tau\nonumber\\
&\leq C(t)\l^{-\f3p_\ast}\Big( \int_0^t\|w[F(\tau,\cdot,\cdot)-F^\v(\tau,\cdot,\cdot)]\|^{\f{b}{b-1}}_{L^2}d\tau\Big)^{\f{2(b-1)}{bp_\ast}} \leq \f{\v}{32}.
\end{align}
For $J_{224132}$, one can obtain, for $\b>3+\max\{0,\g\}$, that 
\begin{equation}\label{3.33}
J_{224132}
\leq C_{N,\v}t\sup_{|v|\leq N,(\tau,y,u)\in((0,T]\times\mathbb{R}^3\times\mathbb{R}^3)}w(u)\Big|F^\v(\tau,y-v\s,u)-F^\v(\tau,y,u)\Big|\rightarrow0,~~\mbox{as}~\s\rightarrow0.
\end{equation}
Combining \eqref{3.24}-\eqref{3.33}, for $\b>6+2\max\{0,\g\}$,  there exists  $\chi_7>0$ depending on $t$ and $\v$ such that if $\s\leq \chi_7$, one has that 
\begin{align}\label{2.60}
J_{22}\leq \f{\v}{32}.
\end{align}
Thus it follows from \eqref{2.11-1}, \eqref{2.11}, \eqref{2.14}, \eqref{2.23} and  \eqref{2.60}, $\b>6+2\max\{0,\g\}$, that 
\begin{align}\label{2.61}
\Big|\int_{\mathbb{R}^3}\Big\{I_2(t+\s,x+\d,v)-I_2(t,x,v)\Big\}dv\Big|\leq \f{\v}{2},~~\mbox{for}~\s+|\d|\leq \min\{\chi_5,\chi_6,\chi_7\}.
\end{align}

Taking $\chi\triangleq \min\{\chi_1,\chi_2,\chi_3,\chi_4,\chi_5,\chi_6,\chi_7\}$, combining \eqref{2.1}, \eqref{3.22-1} and  \eqref{2.61}, for $\b>6+2\max\{0,\g\}$, we have 
\begin{equation}\label{2.62}
\Big|\int_{\mathbb{R}^3}\Big\{F(t+\s,x+\d,v)-F(t,x,v)\Big\}dv\Big|\leq \v,~\mbox{for}~\s+|\d|\leq \chi.
\end{equation}
That is  $\int F(t,x,v) dv$ is continuous function of $(x,t)\in\mathbb{R}^3\times(0,\infty)$. 
In the same way, for $\b>8+2\max\{0,\g\}$, one can prove $\int vF(t,x,v)$ and $\int |v|^2F(t,x,v)dv$ are also continuous functions of $(x,t)\in\mathbb{R}^3\times(0,\infty)$.
Thus the proof of Theorem \ref{thm1.1} is completed.  $\hfill\Box$

\

\noindent{\bf Proof of Theorem \ref{thm1.2}.}
For any fixed  $0<t_1<T+\infty$, let $(x,t)\in\mathbb{R}\times[t_1,T]$.
The uniform continuity of density in $(x,t)\in R^3\times [t_1,T]$ is equivalent to prove that for any small $\v>0$, there exists $\chi>0$ depending on $T,t_1^{-1}$ and $\v$ such that if $|\s|+|\d|\leq \chi$, it holds that
\begin{align}\label{3.1}
\Big|\int_{\mathbb{R}^3}F(t+\sigma,x+\d,v)-F(t,x,v)dv\Big|\leq \v.
\end{align}
 Without loss of generality, we assume $\s\geq0$. As in \eqref{2.4}, one obtains that 
\begin{align}\label{2.4-1}
&K_1 \triangleq\Big|\int_{\mathbb{R}^3}(I_1(t_\s,x+\d,v)-I_1(t,x,v))dv\Big|\nonumber\\
&\leq \int_{\mathbb{R}^3}\Big|F_0(x+\d-vt_\s,v)-F_0(x-vt,v)\Big|dv+\int_{\mathbb{R}^3}F_0(x-vt,v)\Big|\int_t^{t_\s}g(\tau,t_\s,x+\d,v)d\tau\Big|dv\nonumber\\
&~~~~~~~~+\int_0^t \int_{\mathbb{R}^3}F_0(x-vt,v) \Big|g(\tau;t_\s,x+\d,v)-g(\tau;t,x,v)\Big| dvd\tau\nonumber\\
&\triangleq K_{11}+K_{12}+K_{13}.
\end{align}
By the same arguments as in \eqref{2.5}, one has, for $\b>3$, that 
\begin{align}\label{4.1}
K_{11}&\leq \f\v{32}+\sup_{(t,x)\in[t_1,T]\times\mathbb{R}^3}\int_{|v|\leq N}\Big|F_0(x+\d-v(t+\s),v)-F_0(x-vt,v)\Big|dv.
\end{align}
It follows from \eqref{1.5-2} that there exists $\chi_1>0$ depending on $t_1^{-1},T$ and $\v$ such that if $\s+|\d|\leq \chi_1$, one has that 
\begin{equation}\label{4.2}
\sup_{(t,x)\in[t_1,T]\times\mathbb{R}^3}\int_{|v|\leq N}\Big|F_0(x+\d-v(t+\s),v)-F_0(x-vt,v)\Big|dv\leq \f\v{32},
\end{equation}
which, together with \eqref{4.1}, for $\b>3$,  that 
\begin{align}\label{4.3}
K_{11}\leq \f\v{16},~\mbox{for}~\s+|\d|\leq \chi_1.
\end{align}
On the other hand, by the same arguments as in \eqref{2.5-4}-\eqref{3.22}, one can prove that there exists $\chi_2>0$ depending on $t_1^{-1},T$ and $\v$  such that  if $\s+ |\d|\leq \chi_2$, it holds, for $\b>6+2\max\{0,\g\}$,  that 
\begin{align}\label{4.4}
K_{12}+K_{13}\leq \f\v{16}.
\end{align}

As in \eqref{2.11-1}-\eqref{2.61}, one can prove that there exists $\chi_3>0$ depending on $t_1^{-1},T$ and $\v$  such that  if $\s+|\d|\leq \chi_3$, one obtains, for $\b>6+2\max\{0,\g\}$, that 
\begin{align}\label{4.5}
K_{2}\triangleq \Big|\int_{\mathbb{R}^3}I_2(t+\s,x+\d,v)-I_2(t,x,v)dv\Big|\leq \f\v{16}.
\end{align}
Taking $\chi=\min\{\chi_1,\chi_2,\chi_3\}$, then we have, for $\b>6+2\max\{0,\g\}$, that 
\begin{align}\label{4.6}
\Big|\int_{\mathbb{R}^3}F(t+\sigma,x+\d,v)-F(t,x,v)dv\Big|\leq \v~~\mbox{for}~  \s+|\d|\leq \chi.
\end{align}
Thus $\int_{\mathbb{R}^3}F(t,x,v)dv$ is uniformly continuous in $(x,t)\in\mathbb{R}^3\times[t_1,T]$.
In the same way, for $\b>8+2\max\{0,\g\}$, one can prove $\int vF(t,x,v)dv$ and $\int |v|^2F(t,x,v)dv$ are also uniformly continuous in $(x,t)\in\mathbb{R}^3\times[t_1,T]$.
Therefore the proof of Theorem \ref{thm1.2} is completed.  $\hfill\Box$
 
\

\section{Appendix}

We shall prove that the examples given in \eqref{1.8} and \eqref{1.8-1} satisfy the condition \eqref{1.5} and \eqref{1.5-2}. Here we only give the proof for \eqref{1.8} since the case for \eqref{1.8-1} is similar.

Using \eqref{1.9}, it is obvious that $\|w F_0\|_{L^\infty}\leq C<\infty$.
Next we prove the second part of \eqref{1.5}. We only need to prove that for any small $\v>0$, there exists  $\chi>0$ depending  on  $t, \v$ such that if $\s+|\d|\leq \chi$,  it holds that
\begin{align}\label{A4}
\int_{\mathbb{R}^3}|F_0(x+\d-v(t+\s),v)-F_0(x-vt,v)|dv\leq \v.
\end{align}
We note that, for any fixed $t>0$, that 
\begin{align}\label{A1}
&\int_{\mathbb{R}^3}F_0(x+\d-v(t+\s),v)-F_0(x-vt,v)|dv\nonumber\\
&\leq \int_{\mathbb{R}^3}|F_0(x+\d-v(t+\s),v)-F_0(x-v(t+\s),v)|dv\nonumber\\
&~~~~~~~~~~~~~~~~~~~~~~~~~+\int_{\mathbb{R}^3}|F_0(x-v(t+\s),v)-F_0(x-vt,v)|dv\nonumber\\
&\triangleq L_1+L_2.
\end{align}
Without loss of generality, we assume $\s\geq0$. By changing  $y=x-vt_\s$ and \eqref{1.10}, we have
\begin{align}\label{A2}
L_1
&\leq C\Big\{\int_{\mathbb{R}^3}|(\r_0(x+\d-vt_\s)-\r_0(x-vt_\s)|^p+|u_0(x+\d-vt_\s)-u_0(x-vt_\s)|^p\nonumber\\
&~~~~~~~~~~~~~~~~~~~~~~+|\t_0(x+\d-vt_\s)-\t_0(x-vt_\s))|^p dv\Big\}^{\f1p}\nonumber\\
&\leq Ct_\s^{-\f3p}\Big\{\|\r_0(\cdot+\d)-\r_0(\cdot)\|_{L^p}+\|u_0(\cdot+\d)-u_0(\cdot)\|_{L^p}+\|\t_0(\cdot+\d)-\t_0(\cdot))\|_{L^p}\Big\}\nonumber\\
&~~~~~\rightarrow 0, ~~\mbox{as}~\s, \d\rightarrow0,
\end{align}
where we have used the continuity of translations on $L^p$.
Thus there exists $\chi_1>0$ depending on $t>0$ and $\v>0$ such that if $\s+|\d|\leq \chi_1$, it holds that 
\begin{equation}\label{A8}
L_1\leq\f{\v}{2}.
\end{equation}

For $L_2$, it is straightforward to obtain that 
\begin{align}\label{A3}
L_2 &\leq C\int_{\mathbb{R}^3}\Big(|(\r_0(x-vt-v\s)-\r_0(x-vt)|+|u_0(x-vt-v\s)-u_0(x-vt)|\nonumber\\
&~~~~~~~~~~~~~~~~~~~~~~+|\t_0(x-vt-v\s)-\t_0(x-vt))|\Big)e^{-\f{|v|^2}{8\hat{C}}}dv.
\end{align}
 Noting that $\r_0-1\in L^p$ and $0\leq\r_0(x)\leq \hat{C}$, there exists a smooth compact support function $\r_0^\v$ such that
\begin{align}\label{A5}
\|\r_0-\r_0^\v\|_{L^p}\leq \f{\v}{C_2(t)},
\end{align}
where $C_2(t)$ is a large constant. Hence one  obtains that 
\begin{align}\label{A6}
&\int_{\mathbb{R}^3}|\r_0(x-vt-v\s)-\r_0(x-vt)|e^{-\f{|v|^2}{8\hat{C}}}dv\nonumber\\
&\leq \int_{|v|\geq N}|\r_0(x-vt-v\s)-\r_0(x-vt)|e^{-\f{|v|^2}{8\hat{C}}}dv+\int_{|v|\leq N}|\r_0(x-vt-v\s)-\r_0(x-vt)|e^{-\f{|v|^2}{8\hat{C}}}dv\nonumber\\
&\leq Ce^{-\f{N^2}{16\hat{C}}}+ \int_{|v|\leq N}|\r_0(x-vt-v\s)-\r_0^\v(x-vt-v\s)|e^{-\f{|v|^2}{8\hat{C}}}dv\nonumber\\
&~~~~+\int_{|v|\leq N}|\r_0(x-vt)-\r_0^\v(x-vt)|e^{-\f{|v|^2}{8\hat{C}}}dv+\int_{|v|\leq N}|\r_0^\v(x-vt-v\s)-\r_0^\v(x-vt)|e^{-\f{|v|^2}{8\hat{C}}}dv\nonumber\\
&\leq  Ce^{-\f{N^2}{16\hat{C}}}+Ct^{-\f3p}\|\r_0-\r_0^\v\|_{L^p}+\int_{|v|\leq N}|\r_0^\v(x-vt-v\s)-\r_0^\v(x-vt)|e^{-\f{|v|^2}{8\hat{C}}}dv\nonumber\\
&\leq \f{Ct^{-\f3p}}{C_2(t)}\v+Ce^{-\f{N^2}{16\hat{C}}}+C\sup_{y\in\mathbb{R}^3,|v|\leq N}\|\r_0^\v(y+v\s)-\r_0^\v(y)\|_{L^\infty}\leq \f{\v}{16}.
\end{align}
Similarly, we can deal with the remaining terms of \eqref{A3}.  Therefore, there exists  $\chi_2>0$ depending on $t$ and $\v$ such that if $\s\leq \chi_2$, it holds that 
\begin{align}\label{A7}
L_2\leq \f{\v}{4}.
\end{align}
Combining \eqref{A8} and \eqref{A7},  there exists  $\chi=\min\{\chi_1,\chi_2\}$ depending on $t$ and $\v$ such that \eqref{A4} holds if $\s+|\d|\leq \chi$. Therefore the condition \eqref{1.5} holds for the example given  in \eqref{1.8}.

Finally, for any fixed $0<t_1<T<\infty$, from \eqref{A8} and \eqref{A7},  it is not difficult to show  that there exists $\chi>0$ depending on $t_1^{-1},T$ and $\v$ such that if $\s+|\d|\leq \chi$, it holds that 
\begin{equation}\label{A9}
\sup_{(x,t)\in\mathbb{R}^3\times[t_1,T]}\int_{\mathbb{R}^3}|F_0(x+\d-v(t+\s),v)-F_0(x-vt,v)|dv\leq \v.
\end{equation}
Thus the example given in \eqref{1.8}  also satisfies the condition \eqref{1.5-2}.

\

{\bf Acknowledgments.} Feimin Huang is partially supported by by National Center for Mathematics
and Interdisciplinary Sciences, AMSS, CAS and NSFC Grant No.11371349.  Yong Wang is partially supported by NSFC Grant No.  11401565.

\end{document}